\documentclass[12pt]{amsart}

\usepackage{tikz}
\usepackage{tikz-cd}
\usepackage{url,graphicx,tabularx,array}
\usepackage{geometry}
\usepackage{amssymb, amscd, amsfonts}
\usepackage{amsthm}
\usepackage{amsmath, mathtools}
\usepackage{mathrsfs}
\usepackage{wasysym}
\usepackage{latexsym}
\usepackage[all]{xy}
\usepackage{xfrac}
\usepackage{multirow}
\usepackage{listings}
\usepackage{verbatim}
\usepackage{float}
\usepackage{algorithm}
\usepackage{algorithmic}

\usepackage{xcolor}

\usepackage{listings}
\usepackage{color}

\usepackage{lineno}

\usepackage{enumitem}

\definecolor{dkgreen}{rgb}{0,0.6,0}
\definecolor{gray}{rgb}{0.5,0.5,0.5}
\definecolor{mauve}{rgb}{0.58,0,0.82}

\usepackage{array}

\makeatletter
\catcode`! 3
\def\Transpose #1{\romannumeral0\expandafter
                  \Mar@Transpose@a\romannumeral`^^@\Mar@DoOneRow #1\\!\\}

\def\Mar@DoOneRow #1\\{\Mar@DoOneRow@a {}#1&^^@&}%

\def\Mar@DoOneRow@a #1#2&{%
    \if^^@\detokenize{#2}\expandafter\@gobble\fi
    \Mar@DoOneRow@a {#1#2\\}%
}%

\def\Mar@Transpose@a #1#2\\{\ifx!#2\expandafter\Mar@FinishTranspose\fi
    \expandafter\Mar@Transpose@b\romannumeral`^^@\Mar@DoOneRow@a {}#2&^^@&#1}

\def\Mar@Transpose@b #1#2^^@\\{\Mar@Join {}#2^^@!#1}

\def\Mar@Join #1#2\\#3!#4\\%
   {\if^^@\detokenize{#3}\expandafter\Mar@EndJoin\fi
    \Mar@Join {#1#2&#4\\}#3!}%

\def\Mar@EndJoin\Mar@Join #1^^@!^^@\\{\Mar@Transpose@a {#1^^@\\}}

\def\Mar@FinishTranspose
    #1&^^@&#2\\^^@\\{ #2}

\catcode`! 12
\makeatother

\newtheorem{thm}{Theorem}

\newtheorem{lem}[thm]{Lemma}
\newtheorem{prop}[thm]{Proposition}

\theoremstyle{definition}
\newtheorem{ex}[thm]{Example}

\theoremstyle{definition}
\newtheorem{D}[thm]{Definition}

\newcommand{\ds}{\displaystyle}

\let \P \PP

\newcommand{\P}{\mathbb{P}}

\newcommand{\OO}{\mathcal{O}}

\newcommand{\C}{\mathbb{C}}
\newcommand{\R}{\mathbb{R}}

\newcommand{\Z}{\mathbb{Z}}

\newcommand{\g}{\mathfrak{g}}
\newcommand{\fb}{\mathfrak{b}}
\newcommand{\h}{\mathfrak{h}}
\newcommand{\TR}{\textbf{R}}
\newcommand{\BR}{\textbf{\underline{R}}}
\newcommand{\sln}{\mathfrak{sl}_n}
\newcommand{\son}{\mathfrak{so}_n}
\newcommand{\gln}{\mathfrak{gl}_n}
\newcommand{\slt}{\mathfrak{sl}_3}
\newcommand{\slf}{\mathfrak{sl}_4}
\newcommand{\sof}{\mathfrak{so}_4}
\newcommand{\sofive}{\mathfrak{so}_5}
\newcommand{\Mn}{M_{\la n \ra}}

\newcommand{\Seg}{Seg( \P A \times \P B \times \P C)}

\newcommand{\z}{\zeta_6}

\DeclareMathOperator{\codim} {codim}

\newcommand{\la}{\langle}
\newcommand{\ra}{\rangle}

\newcommand{\rk}{\text{rank}}

\newcommand{\Hom}{\text{Hom}}

\usetikzlibrary{arrows,matrix}
\begin{document}

\title{On the Structure Tensor of $\sln$}

\author{Kashif K.~Bari} 
  \email{kashbari@math.tamu.edu} 


\begin{abstract}
The structure tensor of $\sln$, denoted $T_{\sln}$, is the tensor arising from the Lie bracket bilinear operation on the set of traceless $n \times n$ matrices over $\C$. This tensor is intimately related to the well studied matrix multiplication tensor. Studying the structure tensor of $\sln$ may provide further insight into the complexity of matrix multiplication and the ``hay in a haystack'' problem of finding explicit sequences tensors with high rank or border rank. We aim to find new bounds on the rank and border rank of this structure tensor in the case of $\slt$ and $\slf$. {We additionally provide bounds in the case of the lie algebras $\sof$ and $\sofive$.} The lower bounds on the border ranks were obtained via various recent techniques, namely Koszul flattenings, border substitution, and border apolarity. Upper bounds on the rank of $T_{\slt}$ are obtained via numerical methods that allowed us to find an explicit rank decomposition. 

\end{abstract}



\thanks{Author was supported by the Mathematics Department at Texas A{\&}M University and NSF grant AF-1814254.} 

\maketitle

\section{Introduction}

In 1969, Strassen presented a novel algorithm for matrix multiplication of $n \times n$ matrices. {Strassen's algorithm used fewer than the $O(n^3)$ arithmetic operations needed for the standard algorithm.} This led to the question: what is the minimal number of arithmetic operations required to multiply $n \times n$ matrices, or in other words, what is the complexity of matrix multiplication \cite{Str69} \cite{Str83}. An asymptotic version of the problem is to determine the exponent of matrix multiplication, $\omega$, which is the {infimum of all} values such that for all $\epsilon >0$, multiplying $n \times n$ matrices can be performed in $O(n^{\omega+ \epsilon})$ arithmetic operations. Any bilinear operation, including matrix multiplication, may be thought of as tensor in the following way: Let $A,B,$ and $C$ denote vector spaces over $\C$. Given a bilinear map $A^* \times B^* \to C$, the universal property of tensor products induces a linear map $A^* \otimes B^* \to C$. Since $\Hom_{\C}(A^* \otimes B^*, C) \simeq A \otimes B \otimes C$, then we can take our bilinear map to be a tensor in $A \otimes B \otimes C$. Let $\Mn$ denote the matrix multiplication tensor arising from the bilinear operation of multiplying $n \times n$ matrices.  

An important invariant of a tensor is its rank. For a tensor $T \in A \otimes B \otimes C$ the \textit{rank}, denoted $\TR(T)$, is the minimal $r$ such that $T = \sum_{i=1}^r a_i \otimes b_i \otimes c_i$ with $a_i \in A, b_i \in B, c_i \in C$ for $1 \leq i \leq r$. {Given precise $T_i = a_i\otimes b_i\otimes c_i$, then we call $T = \sum_{i=1}^r T_i$ a \textit{rank decomposition} of $T$.} Strassen also showed that the rank of the matrix multiplication tensor is a valid measure of its complexity; in particular, he proved $\omega = \inf \{\tau \in \R \ | \ \TR(\Mn) = O(n^\tau) \}$ \cite{Str69}. 

For a tensor $T \in A \otimes B \otimes C$, the \textit{border rank} of $T$, denoted $\BR(T)$, is another invariant of interest and defined to be the minimal $r$ such that $T = \displaystyle  \lim_{\epsilon \to 0} T_\epsilon$ where for all $\epsilon >0$, $T_\epsilon$ has rank $r$. {Given rank decompositions of $\ds T_\epsilon  = \sum_{i=1}^r T_i(\epsilon)$, we then call $\lim_{\epsilon \to 0} \sum_{i=1}^r T_i(\epsilon)$ a \textit{border rank decomposition} of $T$.} Later, in 1980, it was shown by Bini that the border rank of matrix multiplication is also a valid measure of its complexity by proving that $\omega = \inf \{\tau \in \R \ | \ \BR(\Mn) = O(n^\tau) \}$ \cite{Bini}.

Intimately related to the matrix multiplication tensor is the {structure tensor of the} Lie algebra $\sln$, the set of traceless $n \times n$ matrices over $\C$ equipped with the Lie bracket $[x,y] = xy- yx$. The \textit{structure tensor} of $\sln$ is defined as the tensor arising from the Lie bracket bilinear operation, and we denote it by $T_{\sln}$. {One example of how matrix multiplication is related to $T_{\sln}$ is by closer examination of a skew-symmetric version of the matrix multiplication tensor; consider the tensor arising from the Lie bracket bilinear operation on $\gln$, (which is just $M_n$, but considered as a Lie algebra)} \cite{CHILO}. {Since $\gln = \sln \oplus z$, where $z$ indicates the scalar matrices which are central in $\gln$, $T_{\sln}$ determines the commutator action on all of $\gln$.} While the matrix multiplication tensor has been well studied \cite{CHILO} \cite{CHL}, the structure tensor of $\sln$ has not been studied to the same extent. Currently, the only known non-trivial results are lower bounds on the rank of the structure tensor of $\sln$. {In \cite{GH} it was shown $\TR(T_{\sln}) \geq 2n^2-n-1$}. Studying the structure tensor of $\sln$ may provide some further insight into two central problems in complexity theory.

In complexity theory, it is of interest to find \textit{explicit} objects that behave generically. This type of problem is known as a ``hay in a haystack'' problem. Algebraic geometry tells us that a ``random'' tensor $T$ in $\C^m \otimes \C^m \otimes \C^m$ will have rank/border rank $ \lceil \frac{m^3}{3m-2} \rceil$. {By an \textit{explicit} sequence of tensors,  we will mean a collection} of tensors $T_m \in \C^m \otimes \C^m \otimes \C^m$ such that the coefficients of $T_m$ are computable in polynomial time {in m}. The ``hay in a haystack'' problem for tensors is to find an example of an explicit sequence of tensors of high rank or border rank, {asymptotically in $m$}. Currently, there exists {an explicit} sequence of tensors, $S_m$ {over $\C$}, such that $\TR(S_m) \geq 3m - o(\log(m))$\cite{AFT} and a different {explicit} sequence of tensors {over $\C$}, $T_m$, such that {$\BR(T_m) \geq 2.02m -o(m)$ }\cite{LM}. One should note that the sequence $T_m$ of \cite{LM} has border rank equal to $2m$ when $m = 13$ and has been shown to exceed $2m$ for $m > 364175$. It would be of interest to find sequences of tensors for which the border rank exceeds $2m$ for smaller values of $m$.  

The second problem is Strassen's problem of computing of the complexity of matrix multiplication. The exponent of $\mathfrak{sl}_n$ is defined as  $\displaystyle \omega(  \mathfrak{sl}_n) := \liminf_{n \to \infty} \log_n( \TR(T_{\mathfrak{sl}_n}))$. By Theorem 4.1 from \cite{LY}, the exponent of matrix multiplication is equal to the exponent of $\sln$. Consequently, upper bounds on the rank and even the border rank of $T_{\sln}$ provide upper bounds on $\omega$.

These two problems motivate our study of the border rank of $T_{\sln}$. {We prove new bounds in the case of $\slt$ and $\slf$.}

\begin{thm} \label{sl3_lower}
$\BR(T_{\slt}) \geq 16$
\end{thm}

\begin{thm}\label{sl3_upper}
$\TR(T_{\slt}) \leq  20$ 
\end{thm}

\begin{thm}\label{sl4}
$\BR(T_{\slf}) \geq 28$ 
\end{thm}

{Additionally, we obtain lower bounds on the border rank of the structure tensors of $\sof$ and $\sofive$.} {We note that $\sof \simeq \mathfrak{sl}_2 \times \mathfrak{sl}_2$ and $\sofive \simeq \mathfrak{sp}_4$.}

\begin{thm}\label{so4}
$\BR(T_{\sof}) \geq  9$ 
\end{thm}

\begin{thm}\label{so5}
$\BR(T_{\sofive}) \geq  19$ 
\end{thm}


\section{Preliminaries}

The above definition of $T_{\sln}$ is independent of choice of basis, but we may also write the tensor in terms of bases. Let $\ds \{ a_i \}_{i=1}^{n^2-1}$ be a basis of $\sln$ and $\ds \{ \alpha^i \}_{i=1}^{n^2-1}$ a dual basis.  Recall that $\sln$ has a bilinear operation called the Lie bracket, given by $\ds [a_i, a_j] = a_ia_j - a_ja_i = \sum_{k=1}^{n^2-1} A_{ij}^k a_k$.  The structure tensor of $\sln$ in this basis is $T_{\sln} = \sum_{i,j,k} A_{ij}^k \alpha^i \otimes \alpha^j \otimes a_k \in \sln^* \otimes \sln^* \otimes \sln$. {Let $e_i^j$ denote the $n \times n$ matrix with $1$ in the $(i,j)$th entry. If we take the standard weight basis of $\sln$, namely the matrices $e_i^j$ for $1 \leq i \neq j \leq n$ and $ e_i^i - e_{i+1}^{i+1}$ for $1 \leq i \leq n-1$, then we can express the tensor $T_{\slt}$ as follows.}

\[ 
\begin{bmatrix}
0 & 0 & 0 & 0 & 0 & 0 & 0 & 0 \\
0 & 2 & 0 & 0 & 0 & 0 & 0 & 0 \\
0 & 0 & 1 & 0 & 0 & 0 & 0 & 0 \\
0 & 0 & 0 & -2 & 0 & 0 & 0 & 0 \\
0 & 0 & 0 & 0 & 0 & 0 & 0 & 0 \\
0 & 0 & 0 & 0 & 0 & -1 & 0 & 0 \\
0 & 0 & 0 & 0 & 0 & 0 & -1 & 0 \\
0 & 0 & 0 & 0 & 0 & 0 & 0 & 1
\end{bmatrix} 
,
\begin{bmatrix}
0 & -2 & 0 & 0 & 0 & 0 & 0 & 0 \\
0 & 0 & 0 & 0 & 0 & 0 & 0 & 0 \\
0 & 0 & 0 & 0 & 0 & 0 & 0 & 0 \\
1 & 0 & 0 & 0 & 0 & 0 & 0 & 0 \\
0 & 1 & 0 & 0 & 0 & 0 & 0 & 0 \\
0 & 0 & 1 & 0 & 0 & 0 & 0 & 0 \\
0 & 0 & 0 & 0 & 0 & 0 & 0 & -1 \\
0 & 0 & 0 & 0 & 0 & 0 & 0 & 0
\end{bmatrix} 
\]

\[
\begin{bmatrix}
0 & 0 & -1 & 0 & 0 & 0 & 0 & 0 \\
0 & 0 & 0 & 0 & 0 & 0 & 0 & 0 \\
0 & 0 & 0 & 0 & 0 & 0 & 0 & 0 \\
0 & 0 & 0 & 0 & 0 & -1 & 0 & 0 \\
0 & 0 & -1 & 0 & 0 & 0 & 0 & 0 \\
0 & 0 & 0 & 0 & 0 & 0 & 0 & 0 \\
1 & 0 & 0 & 0 & 1 & 0 & 0 & 0 \\
0 & 1 & 0 & 0 & 0 & 0 & 0 & 0
\end{bmatrix}
,
\begin{bmatrix}
0 & 0 & 0 & 2 & 0 & 0 & 0 & 0 \\
-1 & 0 & 0 & 0 & 0 & 0 & 0 & 0 \\
0 & 0 & 0 & 0 & 0 & 1 & 0 & 0 \\
0 & 0 & 0 & 0 & 0 & 0 & 0 & 0 \\
0 & 0 & 0 & -1 & 0 & 0 & 0 & 0 \\
0 & 0 & 0 & 0 & 0 & 0 & 0 & 0 \\
0 & 0 & 0 & 0 & 0 & 0 & 0 & 0 \\
0 & 0 & 0 & 0 & 0 & 0 & -1 & 0
\end{bmatrix}
\]

\[
\begin{bmatrix}
0 & 0 & 0 & 0 & 0 & 0 & 0 & 0 \\
0 & -1 & 0 & 0 & 0 & 0 & 0 & 0 \\
0 & 0 & 1 & 0 & 0 & 0 & 0 & 0 \\
0 & 0 & 0 & 1 & 0 & 0 & 0 & 0 \\
0 & 0 & 0 & 0 & 0 & 0 & 0 & 0 \\
0 & 0 & 0 & 0 & 0 & 2 & 0 & 0 \\
0 & 0 & 0 & 0 & 0 & 0 & -1 & 0 \\
0 & 0 & 0 & 0 & 0 & 0 & 0 & -2
\end{bmatrix}
,
\begin{bmatrix}
0 & 0 & 0 & 0 & 0 & 1 & 0 & 0 \\
0 & 0 & -1 & 0 & 0 & 0 & 0 & 0 \\
0 & 0 & 0 & 0 & 0 & 0 & 0 & 0 \\
0 & 0 & 0 & 0 & 0 & 0 & 0 & 0 \\
0 & 0 & 0 & 0 & 0 & -2 & 0 & 0 \\
0 & 0 & 0 & 0 & 0 & 0 & 0 & 0 \\
0 & 0 & 0 & 1 & 0 & 0 & 0 & 0 \\
0 & 0 & 0 & 0 & 1 & 0 & 0 & 0
\end{bmatrix}
\]

\[
\begin{bmatrix}
0 & 0 & 0 & 0 & 0 & 0 & 1 & 0 \\
0 & 0 & 0 & 0 & 0 & 0 & 0 & 1 \\
-1 & 0 & 0 & 0 & -1 & 0 & 0 & 0 \\
0 & 0 & 0 & 0 & 0 & 0 & 0 & 0 \\
0 & 0 & 0 & 0 & 0 & 0 & 1 & 0 \\
0 & 0 & 0 & -1 & 0 & 0 & 0 & 0 \\
0 & 0 & 0 & 0 & 0 & 0 & 0 & 0 \\
0 & 0 & 0 & 0 & 0 & 0 & 0 & 0
\end{bmatrix}
,
\begin{bmatrix}
0 & 0 & 0 & 0 & 0 & 0 & 0 & -1 \\
0 & 0 & 0 & 0 & 0 & 0 & 0 & 0 \\
0 & -1 & 0 & 0 & 0 & 0 & 0 & 0 \\
0 & 0 & 0 & 0 & 0 & 0 & 1 & 0 \\
0 & 0 & 0 & 0 & 0 & 0 & 0 & 2 \\
0 & 0 & 0 & 0 & -1 & 0 & 0 & 0 \\
0 & 0 & 0 & 0 & 0 & 0 & 0 & 0 \\
0 & 0 & 0 & 0 & 0 & 0 & 0 & 0
\end{bmatrix} 
\]

We establish some basic definitions of algebraic geometry and representation theory that will be used in our study of the structure tensor of $\sln$.

\subsection{Algebraic Geometry}

The language of algebraic geometry will prove useful in studying our problems. Throughout, let $A,B,C,$ and $V$ be complex vector spaces. {Let $\pi: V \backslash \{0\} \to \P V$ be the projection map from $V$ to its projectivization.} Denote $[v] := \pi(v) \in \P V$ {for nonzero $v \in V$}. {For $v_1, \cdots ,v_k \in V$, let $\la v_1, \cdots, v_k \ra$ denote the projectivization of the linear span of $v_1, \cdots, v_k$.} 

See \cite{Shaf} for the definitions of the Zariski topology, projective variety, and the dimension of a projective variety. {For a projective variety, $X \subset \P V$, let S[X] denote its {homogeneous} coordinate ring and let $I(X) \subset Sym(A^*)$ denote the ideal of $X$.} Given $Y \subset \P V$ a nondegenerate projective variety, we define the $r$th secant variety of $Y$, denoted $\sigma_r(Y) \subset \P V$. 

\begin{D}
The $r$th secant variety of $Y$ is $\sigma_r(Y) = \displaystyle \overline{  \bigcup_{y_i \in Y} \la y_1, \cdots, y_r \ra }$
\end{D}

{The closure operation indicated in Definition 1 is the Zariski closure}. As polynomials are continuous, this closure also contains the Euclidean closure of the set. In this case, we have that the Euclidean closure is in fact equal to the Zariski closure (see Theorem 3.1.6.1 in \cite{GCT} or Theorem 2.3.3 in \cite{Mum}).



%

{Let {$A, B, C$ be vector spaces over $\C$}. The Segre embedding is defined as map $Seg: \P A \times \P B \times \P C \to \P(A \otimes B \otimes C)$ given by $([a] , [b], [c]) \mapsto [a \otimes b \otimes c]$.} Let $X = \Seg$ denote the image of the Segre embedding, which is a projective variety of dimension {$\dim A + \dim B + \dim C - 3$}. Note that $X$ is the space of rank one tensors in $\P(A \otimes B \otimes C)$. Consequently, one can redefine the border rank, {$\BR(T)$, for} $T \in A \otimes B \otimes C$ as the $r$ such that $T \in \sigma_r(X)$ and $T \notin \sigma_{r-1}(X)$. {We make a remark here that $\TR(T) \geq \BR(T)$, trivially, as a tensor rank decomposition can be considered as a border rank decomposition by taking the tensor rank decomposition as a constant sequence.}

For a secant variety $\sigma_r(Y)\subset \P V$ in general, the expected dimension will be $\min \{r \dim Y + (r-1), \dim \P V \}$, since we expect to choose $r$ points from $Y$ and have an additional $r-1$ parameters to span the $r$-plane generated by those points. {The dimension of $\ds \sigma_r(Seg( \P V \times \P V \times \P V))$ is the expected dimension, except when $r=4$ and $\dim V = 3$ \cite{Lic85}.}

%
%

{Another variety we will make use of is the {Grassmannian}, denoted $G(k,V)$. Let $\Lambda^k V$ denote the $k$th exterior power of the vector space $V$. The {Grassmannian} variety is $G(k,V) := \P \{ T \in \Lambda^k V \ | \ \exists v_1, \cdots, v_k \in V \text{ such that } T = v_1 \wedge \cdots \wedge v_k \}$.}




\subsection{Representation Theory}

{Recall that $A,B,C,$ and $V$ are complex vector spaces.} A guiding principle in geometric complexity theory is to use symmetry to reduce the problem of testing a space of tensors to testing particular representatives of families of tensors. To describe the symmetry of our tensor, we use the language of the representation theory of linear algebraic groups and Lie algebras. See \cite{FH} for the definitions of linear algebraic groups, semisimple Lie algebras, representations, orbits of group actions, $G$-modules, irreducibility of a representation/module, Borel subgroups, a maximal torus, and the correspondence between Lie groups and Lie algebras.

The group $GL(A) \times GL(B) \times GL(C)$ acts on $A \otimes B \otimes C$ {by the product of the natural actions of $GL(A)$ on $A$, etc.} {Identify $(\C^*)^{\times 2}$ with the subgroup $\{ (a Id_A, b Id_B, c Id_C) \in GL(A) \times GL(B) \times GL(C) \ | \ abc = 1\}$ of $GL(A) \times GL(B) \times GL(C)$. Note that $(\C^*)^{\times 2}$ acts trivially on $A \otimes B \otimes C$. }

\begin{D}
For a tensor $T \in A \otimes B \otimes C$, define the \textit{symmetry group of} $T$ to be the group $\ds G_T := \{ g \in GL(A) \times GL(B) \times GL(C) / (\C^*)^{\times 2} \ | \ gT = T \}$.

\end{D}

{In the case of $T_{\sln}$, our symmetry group, $G_{T_{\sln}}$, is in fact isomorphic to $SL_n$. For any element $g \in SL_n$, we have the element $g^* \otimes g^* \otimes g$ acts on $\sln^* \otimes \sln^* \otimes \sln$ and leave $T_{\sln}$ invariant. It is always the case that for any automorphism of $\sln$, we will have an automorphism of $\sln^* \otimes \sln^* \otimes \sln$. See \cite{Mir} for a proof that these are all elements of the symmetry group for $T_{\sln}$.}

Let $B_T \subset G_T$ denote a \textit{Borel subgroup}. {In the case of $T_{\sln}$, where our symmetry group {is isomorphic to} $SL_n$, we take $B_{T}$ to be the Borel subgroup of upper triangular matrices of determinant 1. We note that Borel subgroups are not unique, but are all conjugate.} For this Borel subgroup, let $N \subset B_T$ denote the {group} of upper triangular matrices with {diagonal entries equal to 1}, {called the }\textit{{maximal unipotent group}}, and let $\mathbb{T}$ denote the subgroup of diagonal matrices, also called the \textit{maximal torus}. 

\begin{D}
A vector $v \in V^{\otimes k}$ is a \textit{weight vector} if $\mathbb{T} [v] = [v]$. In particular, for $t \in \mathbb{T}$, where $t = \text{diag}\{ t_1, \cdots, t_n \}$, if $ tv = t_1^{p_1} \cdots t_n^{p_n} v$, then $v$ is said to have weight $(p_1, \cdots, p_n) \in \Z^n$. 

Furthermore, call $v$ a highest weight vector if $B_T [v] = [v]$, i.e. if $Nv = v$. 
\end{D}

In our case where $G_T \simeq SL_n$, every irreducible $G_T$-module will have a highest weight line and addtionally will be be uniquely determined by this highest weight line.

Given a symmetry group, $G_T$, we also have a symmetry Lie algebra, denoted $\g_T$, which will be more convenient to work with. Let $\fb_T$ denote the Borel subalgebra and $\h$ denote the Cartan subalgebra, which will be the {Lie algebras} of the Borel subgroup, $B_T$, and maximal torus, $\mathbb{T}$, respectively. In the case of our symmetry Lie algebra, $\sln$, we take $\h$ to be the subalgebra of traceless diagonal matrices. Additionally, for $N \subset B_T$, we have the corresponding Lie subalgebra $\mathfrak{n} \subset \fb_T$, which will consist of the strictly upper triangular elements of $\sln$. We will refer to elements of $\mathfrak{n}$ as raising operators. 

For a Lie algebra representation, a vector is a \textit{weight vector}, if $\h [v] = [v]$. One may regard the weight of a vector $v$ as a linear functional $\lambda \in \h^*$, such that for all $H \in \h$, $Hv = \lambda(H) v$. {Analogous} to the above definition, $[v]$ is a highest weight vector if and only if $\mathfrak{n}  [v] = 0$. Consider the weights $\omega_i \in \h^*$ satisfying $\omega_i( h_j) = \delta_{ij}$ for all $1 \leq i,j \leq n-1$, where $h_j = e_j^j -e_{j+1}^{j+1}  \in \h$. Call these weights $\omega_1, \cdots, \omega_{n-1} \in \mathfrak{h}^*$  the fundamental weights of $\sln$. {It is well known that} the highest weight $\lambda$ of an irreducible representation can be represented as an integral linear combination of fundamental weights. For notational convenience, we denote the irreducible representation with highest weight $\lambda = a_1 \omega_1 + \cdots + a_{n-1} \omega_{n-1}$ by $[a_1 \ \cdots \ a_{n-1}]$ and denote its highest weight vector by $ v_{[a_1 \ \cdots \ a_{n-1}]} $. {We note that this notation is different from the widely used diagram notation of Weyl \cite{FH}}

\begin{D}
 A variety $X \subset \P V$ is $G$\textit{-homogeneous} if it is a closed orbit of some point $x \in \P V$ under the action of some group $G \subset GL(V)$. If $P \subset G$ is the subgroup fixing $x$, write $X = G/P$.
\end{D}


The purpose for introducing the language of homogeneous varieties is to introduce the following Normal Form Lemma. {Recall that $V$ is a complex vector space.} 

\begin{lem}[Normal Form Lemma]
Let $X = G/P \subset \P V$ be a homogeneous variety and $v \in V$ such that $G_v = \{ g \in G \ | \ g[v] = [v] \}$ has a single closed orbit $\OO_{min}$ in $X$. {Then any border rank $r$ decomposition of $v$ may be modified using $G_v$ to a border rank $r$ decomposition $\lim_{\epsilon \to 0}  x_1(\epsilon)+ \cdots+ x_r(\epsilon) $ such that there is a stationary point $x_1(t) \equiv x_1$ (i.e. $x_1$ is independent of $t$) lying in $\OO_{min}$}. 

If, moreover, every orbit of $G_v \cap G_{x_1}$ contains $x_1$ in its closure, we may further assume that for all $j \neq 1$, $\lim_{\epsilon \to 0} x_j(\epsilon) =x_1$.
\end{lem}

See Lemma 3.1 in \cite{NormalForm} for the proof. This Lemma allows one describe the interactions between the different $G$-orbits. It can be thought of as a consequence of the Lie's theorem. {Lie's Theorem states that for a solvable group, $H$, an $H$-module, $W$, and a point $[w] \in \P W$, the closed orbit $\overline{H [w]}$ contains an $H$-fixed point \cite{FH}.}





\section{Methodology}

The most fruitful current techniques for finding lower bounds on the border rank of a tensor are Koszul flattenings, the border substitution method, and border apolarity. We review each of these techniques. 

\subsection{Koszul flattenings}

For $T \in A \otimes B \otimes C$, we may consider it as a linear map $T_B: B^* \to A \otimes C$. {We have analogous maps $T_A, T_B, T_C$, which are called the \textit{coordinate flattenings} of $T$.} Consider the linear map obtained by composing the map $T_B \otimes Id_{\Lambda^p A} : B^* \otimes \Lambda^p A \to \Lambda^{p}A \otimes A \otimes C$ with the {map $\pi \otimes Id_C : \Lambda^{p}A \otimes A \otimes C \to \Lambda^{p+1}A \otimes C$. Note that $\pi \otimes Id_C$ is the tensor product of the exterior multiplication map with the identity on $C$. Denote this composition by $T_A^{\wedge p}$.} Let $rank$ denote the rank of a linear map.

\begin{prop}[Landsberg-Ottaviani]
Let $T \in A \otimes B \otimes C$ and $t = a \otimes b \otimes c \in A \otimes B \otimes C$, then 

\begin{align*}
\displaystyle \BR(T) \geq \frac{ \rk (T_A^{\wedge p})}{\rk ( t_A^{\wedge p})}= \frac{ \rk (T_A^{\wedge p})}{ {\dim A - 1 \choose p }}
\end{align*}
\end{prop}

%
%
%

See \cite{LO} for the proof. One should note that we achieve the best bounds when $\dim A = 2p+1$. Thus, if $\dim A > 2p+1$, we may restrict $T$ to subspaces $A' \subset A$ of dimension $2p+1$, {since border rank is upper semi-continuous with respect to restriction i.e. for a restriction $T'$ of $T$, $\BR(T') \leq \BR(T)$.} Koszul flattenings alone are insufficient to prove $\BR(T_{\sln}) \geq 2(n^2-1)$, as the limit of the method for $T \in \C^m \otimes \C^m \otimes \C^m$ is below $2m-3$(m even) and $2m-5$(m odd). See \cite{LO} for more on this method.


\subsection{Border Substitution}

The only known technique for computing lower bounds on the rank of a tensor is the \textit{substitution method}. A tensor $T$ is \textit{A-concise} if the coordinate flattening map $T_A$ is injective, and define similarly for $B$-concise and $C$-concise. {If a tensor is $A$-concise, $B$-concise, and $C$-concise, then we simply call it \textit{concise}. We remark that $T_{\sln}$ is in fact a concise tensor, since $\sln$ is a simple Lie algebra and the coordinate flattening maps do not send everything to 0.}

%


{The substitution method is a standard technique for obtaining lower bounds on the rank of a tensor \cite{AFT}. This technique can be extended to border rank.}

\begin{prop}[Landsberg-Micha{\l}ek]
Let $T \in A \otimes B \otimes C$ be $A$-concise. Let $\dim A = m$ and let $k < m$. Then 

\begin{align*}
\BR(T) \geq \min_{\ds A' \in G(k,A^*) } \BR( T \big{|}_{A' \otimes B^* \otimes C^*}) + (m-k)
\end{align*}
\end{prop}


{See \cite{LM} for proof. Note that the notation $T \big{|}_{A' \otimes B^* \otimes C^*}$ is a restriction of $T$ when considering $T$ as a trilinear form $T: A^* \otimes B^* \otimes C^* \to \C$. If we let $\tilde{A} = A / (A')^\perp$ then our restricted tensor will be an element of $\tilde{A} \otimes B \otimes C$.}

Also note that in the border substitution proposition, we are minimizing over all elements in the Grassmannian. In practice, applying border substitution uses tensors with large symmetry groups $G_T$. The utility is that one may restrict to looking at representatives of closed $G_T$-orbits in the Grassmannian, rather than by examining all elements of the Grassmannian. One often achieves the best results on the rank of a tensor by using border substitution in conjunction with Koszul flattenings. Naively, {the largest lower bound obtainable by the method, i.e. the limit of the method,} is at most $\dim A + \dim B + \dim C -3$, however, the limit is in fact slightly less. For tensors $T \in \C^m \otimes \C^m \otimes \C^m$, the limit of the method is $3m - 3\sqrt{3m + \frac{9}{4}} + \frac{9}{2}$. See \cite{NormalForm} for a proof of this and more on this method. The best lower bound achieved on the border rank that is mentioned in the introduction is achieved using this method \cite{LM}.

\subsection{Border Apolarity}
 
{In order to establish larger lower bounds on $\BR(T_{\slt})$ than can be achieved by Koszul flattenings and border substitution for $T_{\slt}$, we will use {border apolarity}, as developed in \cite{BB} and \cite{CHL}.}
 
Suppose $T$ has a border rank $r$ decomposition, $ T = \lim_{\epsilon \to 0} T_{\epsilon}$, where $ T_\epsilon  = \sum_{i=1}^r T_i(\epsilon)$. {If the rank summands $T_i(\epsilon)$ are in general position in $A \otimes B \otimes C$, then we may identify the border rank decomposition with a curve $E_\epsilon$ in the Grassmannian variety $G(r, A\otimes B \otimes C)$, by taking the exterior product of the $T_i(\epsilon)$, i.e. $E_\epsilon= [ T_1(\epsilon) \wedge \cdots \wedge T_r(\epsilon)]$.}

{Now we define a $\Z^3$-grading on ideals of subsets of $\P A \times \P B \times \P C$ (i.e. the Segre variety) from the natural $\Z^3$-grading of $Sym (A \oplus B \oplus C)^*$. Let $Irrel := \{0 \oplus B \oplus C \} \cup \{ A \oplus 0 \oplus C \} \cup \{ A \oplus B \oplus 0 \} \subset A \oplus B \oplus C$. Since $\P A \times \P B \times \P C \simeq (A \oplus B \oplus C \backslash Irrel) / (\C^*)^{\times 3}$, then we may consider the quotient map $q: (A \oplus B \oplus C) \backslash Irrel \to \P A \times \P B \times \P C$, which will be invariant under the action of $(\C^*)^{\times 3}$. Therefore, for a set $Z \subset \P A \times \P B \times \P C)$, the ideal of this set $I(Z) = I(q^{-1}(Z)) \subset Sym(A \oplus B \oplus C)^*$ will have a $\Z^3$ grading. In particular, for a single point $([a], [b], [c]) \in \P A \times \P A \times \P C$, corresponding to a rank one tensor ($[a\otimes b \otimes c] \in \P ( A \otimes B \otimes C)$), we are considering the ideal in $Sym(A \oplus B \oplus C)^*$ of polynomials vanishing along the lines $a,b,$ and $c$.} 

Let $I_\epsilon$ denote the $\Z^3$-graded ideal of the {set of the} r points $[ T_i( \epsilon)]$, i.e. $I_{ijk,\epsilon}\subset S^iA^* \otimes S^j B^* \otimes S^k C^*$. {Since the $r$ points are in general position, then $\codim I_{ijk,\epsilon} =r$.} Define $\ds I_{ijk} := \lim_{\epsilon \to 0} I_{ijk,\epsilon}$ as the limit of points in the Grassmannian $G(\dim(S^iA^* \otimes S^jB^* \otimes S^kC^*) - r,\dim(S^iA^* \otimes S^jB^* \otimes S^kC^*) )$. $I_{ijk}$ will exist, since the Grassmannian is compact, however, the {resulting ideal $I$ may not be saturated}. See \cite{BB} for further discussion on this. 

{Recall that a tensor $T$ is concise if all the coordinate flattening maps $T_A, T_B,T_C$ are injective. For a subspace $U \subset V$, define $U^\perp := \{ \alpha \in V^* \ | \ \alpha(u) = 0 \ \forall u \in U \}$. In a nutshell, border apolarity gives us some necessary conditions on the possible limiting ideals, $I$, that can arise from a border rank decomposition.}

\begin{thm}(Weak Border Apolarity)
Let $X = \P A \times \P B \times \P C$ and $S[X]$ be its coordinate ring. Suppose a tensor $T$ has $\BR(T) \leq r$. Then there exists a (multi)homogeneous ideal $I \subset S[X]$ such that

\begin{itemize}

\item $I \subset Ann(T) $ 

\item For each multidegree $ijk$, the $ijk$th graded piece, $I_{ijk}$, of $I$ has $\codim I_{ijk} = \min\{r, \dim S[X]_{ijk}\}$ 
\end{itemize}

In addition, if $G_T$ is a group acting on $X$ and preserving $T$, then there exists an $I$ as above which in addition is invariant under a Borel subgroup of $G_T$.

\end{thm} 

{In \cite{BB}, see Theorems 3.15 (Border Apolarity) for proof of the first part and see Theorem 4.3 (Fixed Ideal Theorem) for a proof of the second part. Theorem 3.15 and 4.3 of \cite{BB} are not stated here as they are stated in greater generality than we require using the language of schemes. We remark that the Weak Border Apolarity Theorem provides sufficiency that if a border rank $r$ decomposition exists, then there will exist another border rank $r$ decomposition satisfying the given conditions.}

{We note that the second condition says that we may in fact take the $r$ points of the border rank decomposition to be in general position, and so our initial supposition that the $r$ points are in general position is justified. Lie's Theorem and the Normal Form Lemma allow us to take $I_{111}$ to be $B_T$-fixed. The Fixed Ideal Theorem of \cite{BB} uses the same reasoning to generalize this to prove $B_T$-invariance for all multigraded components $I_{ijk}$, not just a finite number of multigraded components.}

{In \cite{CHL}, using Weak Border Apolarity Theorem, they assert that for $T$ a concise tensor with a border rank $r$ decomposition, there will exist an ideal $I$ satisfying the following:}

\begin{enumerate}
	\item $I_{ijk}$ is $B_T$-stable ($I_{ijk}$ is a Borel fixed weight space)
	\item $I \subset Ann(T) $ i.e. $I_{110} \subset T(C^*)^{\perp} \subset A^* \otimes B^*$, etc. and $I_{111} \subset T^{\perp}\subset A^* \otimes B^* \otimes C$
	\item For all $i,j,k$ such that $ i+j+k >1$, then $\codim I_{ijk} = r$, i.e. the condition that we may take the $r$ points of the border rank decomposition to be in general position.
	\item Since $I$ is an ideal, the image of the multiplication map $I_{i-1,j,k} \otimes A^* \oplus I_{i,j-1,k} \otimes B^* \oplus I_{i,j,k-1} \otimes C^* \to S^iA^* \otimes S^jB^* \otimes S^kC^*$ is contained in $I_{ijk}$
\end{enumerate}

{We note that the last condition is simply the condition that $I$ is an ideal where the multiplication respects the grading. The border apolarity algorithm presented in \cite{CHL} makes use of these conditions and attempts to iteratively construct all possible ideals, $I$, in each multidegree. In particular, if at any multidegree $ijk$, there does not exist an $I_{ijk}$ satisfying the above, then we may conclude that $\BR(T) > r$. We remark that the candidate ideals, $I$, do not necessarily correspond to an actual border rank decomposition of the tensor. We describe the algorithm of \cite{CHL} precisely below:}

\begin{algorithm}[H]
\caption{Border Apolarity Algorithm of Conner, Harper, Landsberg}
Input: $T \in A \otimes B \otimes C$, $r$ 

Output:  Candidate ideals or $\BR(T) > r$

\begin{enumerate}
 \item For each $B_T$-fixed space $F_{110}$ of codim $r- \dim C$ in $T(C^*)^\perp$ (i.e. codim $r$ in $A^* \otimes B^*$) compute ranks of maps $F_{110} \otimes A^* \to S^2 A^* \otimes B^*$ and $F_{110} \otimes B^* \to A^* \otimes S^2B^*$

If both have images of codim at least $r$, then $F_{110}$ is possible $I_{110}$. These are called the (210) and (120) tests.

 \item Perform analogously for possible $I_{101} \subset T(B^*)^\perp$ and $I_{011} \subset T(A^*)^\perp$ for candidate $F_{101}$ and $F_{011}$

 \item For each triple $F_{110}, F_{101}, F_{011}$, compute rank of map $(F_{110} \otimes C^*) \oplus (F_{101} \otimes B^*) \oplus (F_{011} \otimes A^*) \to A^* \otimes B^* \otimes C^*$

If codim of image is at least $r$, then have a candidate triple. $F_{111}$ is candidate for $I_{111}$ if it is codim $r$, it is contained in $T^\perp$ and contains image of above map.

\item  Analogous higher degree tests



\item  If at any point there are no such candidates $\BR(T) >r$, otherwise stabilization of candidate ideals will occur at worst multi-degree $(r,r,r)$
 \end{enumerate}
\end{algorithm}

{The condition that $I_{ijk}$ is $B_T$-fixed allows us to greatly reduce the search for possible candidate ideals. The $B_T$-fixed spaces are easier to list than trying to list all possible $I_{ijk}$. This condition allows the algorithm of \cite{CHL} to be feasible for tensors with large symmetry groups. Then, using the assumption that our points are in general position, one has rank conditions on the multiplication maps, as the images must have codim at least $r$.}

\subsubsection{{Implementation of Border Apolarity for $T_{\slt}$}}

We show how to implement the algorithm of \cite{CHL} for $T_{\sln}$ by describing how to compute all possible $B_T$-fixed $I_{110}$. Additionally, we can leverage the skew-symmetry of $T_{\sln}$ to reduce the amount of computation involved for determining potential $I_{110}, I_{101}, I_{011}$.

Let $T = T_{\mathfrak{sl}_3} \in \slt^* \otimes \slt^* \otimes \slt = A \otimes B \otimes C$. The first and third condition from border apolarity tells us to compute all $B_T$-fixed weight subspaces $F_{110} \subset A^* \otimes B^*$ of codimension $r$; however, since $T$ is concise with $\dim T(C^*) = \dim \slt = 8$ and $F_{110} \subset T(C^*)^\perp$ by the second condition of border apolarity, we compute all $B_T$-fixed weight spaces $F_{110} \subset T(C^*)^\perp$ of codimension $r-8$. 

{Standard computational methods, see \cite{FH}, yield that the irreducible decomposition of $A^* \otimes B^*$ as $\sln$-modules is as follows.}

\begin{align*}
&A^* \otimes B^* = \sln \otimes \sln \\
&\simeq  [2 \ 0 \ \cdots \ 0 \ 2] \oplus [2 \ 0 \cdots 0 \ 1 \ 0] \oplus [0 \ 1 \ 0 \cdots \ 0 \ 2]\oplus [0 \ 1 \ 0 \cdots \ 0 \ 1 \ 0] \oplus 2[1 \ 0 \ \cdots 0 \ 1] \oplus [0 \ \cdots \ 0]
\end{align*} 

 In particular, for $\slt$, we have the following decomposition into $\slt$-modules.
 
 \begin{align*}
 \slt \otimes \slt \simeq [2 \ 2] \oplus [3 \ 0] \oplus [0 \ 3] \oplus 2[1 \ 1] \oplus [0 \ 0]
 \end{align*}
 
\begin{figure*}[ht]
\centering
\includegraphics[scale=0.75]{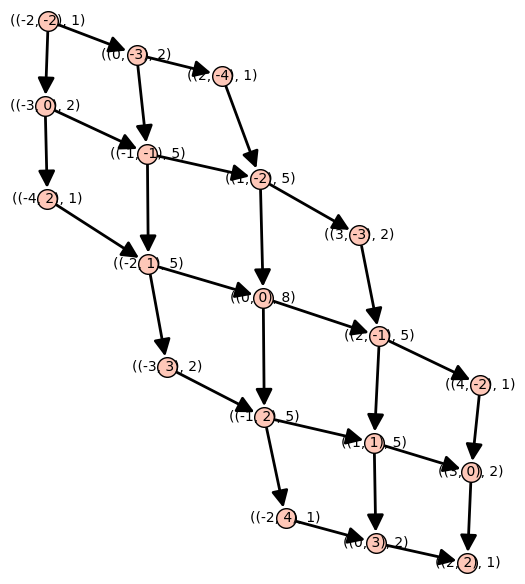}
\caption{Weight Decomposition of $T(C^*)^\perp$}
\label{fig:poset}
\end{figure*}

   Using this decomposition and the conciseness of $T$, we have the $\slt$-module decomposition $T(C^*)^\perp \simeq [2 \ 2] \oplus [3 \ 0] \oplus [0 \ 3] \oplus [1 \ 1] \oplus [0 \ 0]$ (Note that $\dim T(C^*)^\perp = 56$). Using this $\slt$-module decomposition, we can obtain a decomposition of $T(C^*)^\perp$ into weight spaces (decomposition into $\h$-modules) by {combining the weight space decompositions} of each $\slt$-module into one poset. The result is the Figure~\ref{fig:poset}. Each node represents a weight space of weight $\lambda$ labeled by $(\lambda,$ dimension of weight space in $T(C^*)^\perp)$. {Also note that arrows go from lower weights to higher weights, so the highest weight occuring in $T(C^*)^\perp$ is $[2 \ 2]$. The weight space $[3 \ 0]$ is of dimension 2, where one of the basis elements for this weight space is the highest weight vector of the $\slt$-module $[3 \ 0]$ and the other basis element for this weight space is the vector arising from lowering the highest weight vector of the $\slt$-module $[2 \ 2]$. }

The utility of this decomposition is that we can generate all $B_T$-fixed weight subspaces $F_{110}$ by taking a collection of weight vectors $v_i$ from the poset such that $v_1 \wedge v_2 \wedge \cdots \wedge v_{56-(r-8)}$ is a highest weight vector in $G(56-(r-8),T(C^*)^\perp)$ and consequently is closed under the raising operators. One should note that if $v_i$ comes from a weight space of dimension greater than 1 then one needs to include linear combinations of basis vectors of that weight space. 

\begin{ex}
{We show a few small examples of possible $F_{110}$ $B_T$-fixed subspaces to provide intuition of how these spaces are computed.} 

For $r = 63$, we have that $F_{110}$ will be of the form $v_1$. Since it must be closed under raising operators, it will necessarily be a highest weight vector. Therefore, our choices will be $v_1 = v_\lambda$ where of $v_\lambda$ is a highest weight vector of weight $\lambda = [2 \ 2]$, $[3 \ 0]$, $[0 \ 3]$, $[1 \ 1]$, or $[0 \ 0]$.

For $r = 62$, $F_{110}$ will be of the form $v_1 \wedge v_2$. Necessarily, we must have that $v_1$ must be a highest weight vector. The second vector $v_2$ may either be another highest weight vector, a weight vector that can be raised to $v_1$, or a {linear combination of the two previous cases if they are vectors of the same weight}. For example, if we take $v_1 = v_{[2 \ 2]}$ to be the highest weight vector of $[ 2 \ 2]$. Let $v_{[3 \ 0]}$ and $u_{[3 \ 0]}$ be a weight basis for $[3 \ 0]$ with $v$ being a highest weight vector. Assume similarly for $[0 \ 3]$ weight space. {The possible choices for $v_2$ are weight vectors of the following types:} $v_{[3 \ 0]}, sv_{[3 \ 0]}+tu_{[3 \ 0]}, v_{[0 \ 3]}, sv_{[0 \ 3]}+tu_{[0 \ 3]}, v_{[1 \ 1]}, v_{[0 \ 0]}$ with $s,t$ parameters. Applying all possible raising operators to $v_1 \wedge v_2$ where $v_2$ has parameters will provide equations for what values of $s,t$ give us a highest weight vector.
\end{ex}

{For smaller values of $r$, the number of possible Borel fixed spaces is much larger and more difficult to list by hand without the aide of a computer.} The computationally difficult step in this algorithm lies in computing the ranks of the multiplication maps such as $F_{110} \otimes A^* \to S^2A^* \otimes B^*$. In some cases there are many parameters which arise from choosing weight vectors from high dimensional weight spaces, such as the $[0 \ 0]$ weight space in Figure~\ref{fig:poset}. Recall that a linear map has rank at most $k$ is the $k+1$ minors all vanish. In order to determine whether the multiplication map has image codimension $r$, we look at the appropriate minors of this linear map. When there are no parameters involved, this is a simple linear algebra calculation. However, in some cases, the entries of the multiplication map are linear polynomials in the parameters coming from choosing a linear combination of weight vectors. In order to determine whether the multiplication map has image of codimension $r$, one needs to look at the ideal of these minors, as well as some polynomial equations in the parameters that are needed for the space to be Borel fixed. One must do a Groebner Basis computation on this ideal to determine whether all the minors vanish or not. This can become an unfeasible computation if there are too many parameters and/or equations.

\section{New Bounds}

{It is known that $\BR( T_{\mathfrak{sl}_2} ) = 5$ \cite{GH}, so we aim to find bounds on $\BR( T_{\sln} ) $ for $n = 3$ and $4$ using the above techniques.}{ Additionally, we provide new lower bounds on $\BR(T_{\son})$ for $n =4$ and $5$.}

\subsection{Koszul flattenings}

In the case of $T_{\slt}$, we achieve the best results when $p=3$ and we restrict to a generic $7$ dimensional subspace of $\slt$, since $\dim \slt = 8$. The best bound achieved is $\BR(T_{\slt}) \geq 14$ (See Table 2).  


\begin{table}[H]
\caption{$T_{\mathfrak{sl}_3}$ Results} 
\centering 
\begin{tabular}{c c c c} 
\hline\hline 
p & Dimensions of Linear Map & Dimension of Kernel & Koszul Bound \\ [0.5ex] 
\hline 
1 & (64,224) & 0 & 10 \\ 
2 & (224,448) & 1 &  11 \\ 
3 & (448,560) & 8 & 13 \\ 
\hline 
\end{tabular}
\end{table}
%
%

\begin{table}[H]
\caption{$T_{\mathfrak{sl}_3}$ Restriction to a generic subspace of dim $k$ Results} 
\centering 
\begin{tabular}{c c c c c} 
\hline\hline 
p & k & Dimensions of Linear Map & Dimension of Kernel & Koszul Bound \\ [0.5ex] 
\hline 
1 & 3 &  (24,24) & 0 & 12 \\
2 & 5 & (80,80) & 4 & 13 \\ 
3 & 7 & (280,280) & 7 & 14 \\ 
\hline 
\end{tabular}
\end{table}

%

In the case of $T_{\slf}$, we achieve the best lower bound of $27$ when $p=4$ or $5$ while restricting to a subspace (See Table 4).

\begin{table}[H]
\caption{$T_{\mathfrak{sl}_4}$ Results} 
\centering 
\begin{tabular}{c c c c} 
\hline\hline 
p & Dimensions of Linear Map & Dimension of Kernel & Koszul Bound \\ [0.5ex] 
\hline 
1 & (225,1575) & 0 & 17\\ 
2 & (1575,6825) & 1 &  18 \\
3 & (6825,20475) & 15 &  19\\ 
4 & (20475,45045) & 106 &  21 \\ 
5 & (45045,75075) & 470 &  23\\
6 & (75075,96525) & 2680 &  25\\  
7 & (96525,96525) & 11039 &  25\\ [1ex] 
\hline 
\end{tabular}
\end{table}
%
%

\begin{table}[H]
\caption{$T_{\mathfrak{sl}_4}$ Restriction to a generic subspace of dim $k$ Results } 
\centering 
\begin{tabular}{c c c c c} 
\hline\hline 
p & k & Dimensions of Linear Map & Dimension of Kernel & Koszul Bound \\ [0.5ex] 
\hline 
1 & 3 & (45,45) & 0 & 23\\ 
2 & 5 & (150,150) & 2 &  25 \\ 
3 & 7 & (525,525) & 7 &  26\\  
4 & 9 & (1890,1890) & 38 &  27 \\ 
5 & 11 & (6930,6930) & 176 &  27\\  
6 & 13 & (25740,25740) &  2254 &  26 \\ [1ex] 
\hline 
\end{tabular}
\end{table}

As stated above, Koszul flattenings alone are insufficient to obtain border rank lower bounds exceeding $2m$, i.e. Koszul flattenings will not prove $\BR( T_{\slt}) \geq 16$ and $\BR(T_{\slf}) \geq 30$.

{In addition to lower bounds on these tensor, we computed lower bounds on $T_{\sof}$ and $T_{\sofive}$.}

\begin{table}[H]
\caption{$T_{\mathfrak{so}_4}$ Results} 
\centering 
\begin{tabular}{c c c c} 
\hline\hline 
p & Dimensions of Linear Map & Rank & Koszul Bound \\ [0.5ex] 
\hline 
1 & (90,36) & 36 & 8\\ 
2 & (120,90) & 78 &  8 \\ 
3 & (90,120) & 78 &  8 \\ [1ex]
\hline 
\end{tabular}
\end{table}

\begin{table}[H]
\caption{$T_{\mathfrak{so}_4}$ Restriction to a generic subspace of dim $k$ Results } 
\centering 
\begin{tabular}{c c c c c} 
\hline\hline 
p & k & Dimensions of Linear Map & Rank & Koszul Bound \\ [0.5ex] 
\hline 
1 & 3 & (18,18) & 18 & 9\\ 
2 & 5 & (60,60) & 48 &  8 \\ [1ex] 
\hline 
\end{tabular}
\end{table}

{From \cite{Mir}, it is shown that for the Lie algebras $\mathfrak{g}_n = \mathfrak{sl}_2^ {\times n}$ that $\TR(T_{\mathfrak{g}_n}) = 5n$, so  $\TR(T_{\sof}) = 10$. Therefore, $\BR(T_{\sof})$ will either be 9 or 10.}

\begin{table}[H]
\caption{$T_{\mathfrak{so}_5}$ Results} 
\centering 
\begin{tabular}{c c c c} 
\hline\hline 
p & Dimensions of Linear Map & Rank & Koszul Bound \\ [0.5ex] 
\hline 
1 & (450,100) & 100 & 12\\ 
2 & (1200,450) & 449 &  13 \\
3 & (2100,1200) & 1190 &  15\\ 
4 & (2520,2100) & 1971 &  16 \\ 
5 & (2100,2520) & 1971 &  16 \\ [1ex] 
\hline 
\end{tabular}
\end{table}

\begin{table}[H]
\caption{$T_{\mathfrak{so}_5}$ Restriction to a generic subspace of dim $k$ Results } 
\centering 
\begin{tabular}{c c c c c} 
\hline\hline 
p & k & Dimensions of Linear Map & Rank & Koszul Bound \\ [0.5ex] 
\hline 
1 & 3 & (30,30) & 30 & 15\\ 
2 & 5 & (100,100) & 98 &  17 \\ 
3 & 7 & (350,350) & 343 &  18\\  
4 & 9 & (1260,1260) & 1136 &  17 \\  [1ex] 
\hline 
\end{tabular}
\end{table}

\subsection{Border Substitution}

{For $T_{\sln} \in \sln^* \otimes \sln^* \otimes \sln$, we may identify the space $\sln^*$ with $\sln$ (by sending an element to its negative transpose). Therefore, we may identify $T_{\sln}$ as an element of $\sln \otimes \sln \otimes \sln$. As a first step in applying border substitution, we restrict $T_{\sln} \in A \otimes B \otimes C$ in the $A$ tensor factor.} Since we may restrict to looking at representatives of closed $G_{T_{\sln}}$-orbits, then the only planes we need to check are the highest weight planes in $G(k, \sln)$. In order to compute the border rank of the restricted tensor, we use Koszul flattenings on the restricted tensor. Once again, let $v_\lambda$ denote the unique weight vector in weight space $\lambda$. 

For $T_{\slt}$, border substitution did not generate a better lower bound than the Koszul flattenings. However, we were able to obtain a better lower bound for $T_{\slf}$. Let $A'$, as in Proposition 15, be $\tilde{A}^\perp$ where we take $\tilde{A}$ to be a space of dimension $m-k$. If we restrict our tensor by a one dimensional subspace, then the only choice for $\tilde{A}$ will be the space spanned by $v_{[1 \ 0 \ 1]}$, which is the highest weight vector of $\sln$.

\begin{table}[H]
\caption{$T_{\slf}$ with Restriction $\tilde{A} = v_{[1 \ 0 \ 1]}$} 
\centering 
\begin{tabular}{c c c c c} 
\hline\hline 
p & k & Dimensions of Linear Map & Dimension of Kernel & Koszul Bound \\ [0.5ex]
\hline 
1 & 3 & (45,45) & 0 & 23\\ 
2 & 5 & (150,150) & 2 &  25 \\ 
3 & 7 & (525,525) & 7 &  26\\  
4 & 9 & (1890,1890) & 38 &  27 \\ 
5 & 11 & (6930,6930) & 248 &  27\\  
6 & 13 & (25740,25740) &  2254 &  26 \\ [1ex] 
\hline 
\end{tabular}
\end{table}

Restricting by a two dimensional subspace, we have one choice for $\tilde{A}$ up to symmetry in the weight space decomposition for $\slf$, namely $v_{[1 \ 0 \ 1]}  \wedge v_{[-1 \ 1 \ 1]}$.

\begin{table}[H]
\caption{$T_{\slf}$ with Restriction  $\tilde{A} = v_{[1 \ 0 \ 1]}  \wedge v_{[-1 \ 1 \ 1]}$} 
\centering 
\begin{tabular}{c c c c c} 
\hline\hline 
p & k & Dimensions of Linear Map & Dimension of Kernel & Koszul Bound \\ [0.5ex]
\hline 
1 & 3 & (45,45) & 0 & 23\\ 
2 & 5 & (150,150) & 2 &  25 \\ 
3 & 7 & (525,525) & 7 &  26\\  
4 & 9 & (1890,1890) & 78  &  26 \\ 
5 & 11 & (6930,6930) &  498 &  26 \\ [1ex] 
\hline 
\end{tabular}
\end{table}

Restricting by a three dimensional subspace, we have three choices for $\tilde{A}$ up to symmetry. $\tilde{A}$ may be $v_{[1 \ 0 \ 1]}  \wedge v_{[1 \ 1 \ -1]} \wedge v_{[-1 \ 1 \ 1 ]}$, $v_{[1 \ 0 \ 1]}  \wedge v_{[1 \ 1 \ -1]} \wedge v_{[2 \ -1 \ 0 ]}$, or $v_{[1 \ 0 \ 1]}  \wedge v_{[1 \ 1 \ -1]} \wedge v_{[-1 \ 2 \ -1 ]}$.


\begin{table}[H]
\caption{$T_{\slf}$ with Restriction $\tilde{A} =  v_{[1 \ 0 \ 1]}  \wedge v_{[1 \ 1 \ -1]} \wedge v_{[-1 \ 1 \ 1 ]}$} 
\centering 
\begin{tabular}{c c c c c} 
\hline\hline 
p & k & Dimensions of Linear Map & Dimension of Kernel & Koszul Bound \\ [0.5ex]
\hline 
1 & 3 & (45,45) & 0 & 23\\ 
2 & 5 & (150,150) & 2 &  25 \\ 
3 & 7 & (525,525) & 31 &  25\\  
4 & 9 & (1890,1890) & 168  &  25 \\ 
5 & 11 & (6930,6930) &  755 &  25 \\ [1ex] 
\hline 
\end{tabular}
\end{table}


\begin{table}[H]
\caption{$T_{\slf}$ with Restriction $\tilde{A} =  v_{[1 \ 0 \ 1]}  \wedge v_{[1 \ 1 \ -1]} \wedge v_{[2 \ -1 \ 0 ]}$} 
\centering 
\begin{tabular}{c c c c c} 
\hline\hline 
p & k & Dimensions of Linear Map & Dimension of Kernel & Koszul Bound \\ [0.5ex]
\hline 
1 & 3 & (45,45) & 0 & 23\\ 
2 & 5 & (150,150) & 2 &  25 \\ 
3 & 7 & (525,525) & 7 &  25\\  
4 & 9 & (1890,1890) & 72  &  25 \\ 
5 & 11 & (6930,6930) &  498 &  25 \\ [1ex] 
\hline 
\end{tabular}
\end{table}


\begin{table}[H]
\caption{$T_{\slf}$ with Restriction $\tilde{A} = v_{[1 \ 0 \ 1]}  \wedge v_{[1 \ 1 \ -1]} \wedge v_{[-1 \ 2 \ -1 ]}$} 
\centering 
\begin{tabular}{c c c c c} 
\hline\hline 
p & k & Dimensions of Linear Map & Dimension of Kernel & Koszul Bound \\ [0.5ex]
\hline 
1 & 3 & (45,45) & 3 & 21\\ 
2 & 5 & (150,150) & 14 &  23 \\ 
3 & 7 & (525,525) & 42 &  25\\  
4 & 9 & (1890,1890) & 254  &  24 \\ 
5 & 11 & (6930,6930) &  1072 &  24 \\ [1ex] 
\hline 
\end{tabular}
\end{table}

The best bound we obtain is $\BR(T_{\slf} \big{|}_{(v_{[1 \ 0 \ 1]})^\perp \otimes \sln \otimes \sln}) \geq 27$ (See Table 9). By Proposition 15, this proves Theorem \ref{sl4}, $\BR(T_{\slf}) \geq 28$. After restricting in the $A$ tensor factor, we cannot restrict by the highest weight vector of $\sln$ in the $B$ or $C$ factor as the symmetry group of the restricted tensor will have a different symmetry group.

{Additionally, border substitution gives a better lower bound on $\BR(T_{\sofive})$. Restricting by a two dimension subspace (See Table 14) and using Proposition 15 proves Theorem \ref{so5}.}

\begin{table}[H]
\caption{$T_{\sofive}$ with Restriction $\tilde{A} = v_{[0 \ 1]}  \wedge v_{[1 \ 0]}$} 
\centering 
\begin{tabular}{c c c c c} 
\hline\hline 
p & k & Dimensions of Linear Map & Dimension of Kernel & Koszul Bound \\ [0.5ex]
\hline 
1 & 3 & (30,30) & 0 & 15\\ 
2 & 5 & (100,100) & 4 &  16 \\ 
3 & 7 & (350,350) & 28 &  17\\  
4 & 9 & (1260,1260) & 194  &  16 \\ [1ex]
\hline 
\end{tabular}
\end{table}

\subsection{Border Apolarity}

We use border apolarity to disprove that $T_{\slt}$  has rank $r =15$. We first compute candidate $F_{110}$ spaces which passed the $(210)$-test. There were a total of 5 candidate $F_{110}$ subspaces out of a total of more than 1245 possible $F_{110}$ spaces. The candidate $F_{110}$ spaces came in three types of weight space decompositions:

\begin{table}[H]
\caption{Two candidate $F_{110}$ planes have the following weight decomposition} 
\centering 
\begin{tabular}{c c c} 
\hline\hline 
Weight & Dimension of weight space in $A^* \otimes B^*$ & Dimension of weight space in $F_{110}$  \\ [0.5ex]
\hline 
$[2, 2]$ & 1 & 1 \\ 
$[3, 0]$ & 2 & 2 \\ 
$[4 , -2]$ & 1 & 1 \\  
$[0 , 3]$ & 2 & 2 \\ 
$[1 , 1]$ & 5 & 5 \\ 
$[2 , -1]$ & 5 & 5 \\ 
$[3 , -3]$ & 2 & 2 \\  
$[-2 , 4]$ & 1 & 1 \\ 
$[-1 , 2]$ & 5 & 5 \\ 
$[0 , 0]$ & 8 & 8 \\ 
$[1 , -2]$ & 5 & 5 \\  
$[2 , -4]$ & 1 & 1 \\ 
$[-3 , 3]$ & 2 & 2 \\ 
$[-2 , 1]$ & 5 & 4 \\ 
$[-1 , -1]$ & 5 & 4 \\  
$[0 , -3]$ & 2 & 1 \\ [1ex] 
\hline 
\end{tabular}
\end{table}

\begin{table}[H]
\caption{One candidate $F_{110}$ plane has the following weight decomposition} 
\centering 
\begin{tabular}{c c c} 
\hline\hline 
Weight & Dimension of weight space in $A^* \otimes B^*$ & Dimension of weight space  in $F_{110}$ \\ [0.5ex]
\hline 
$[2 , 2]$ & 1 & 1 \\ 
$[3 , 0]$ & 2 & 2 \\ 
$[4 , -2]$ & 1 & 1 \\  
$[0 , 3]$ & 2 & 2 \\ 
$[1 , 1]$ & 5 & 5 \\ 
$[2 , -1]$ & 5 & 5 \\ 
$[3 , -3]$ & 2 & 2 \\  
$[-2 , 4]$ & 1 & 1 \\ 
$[-1 , 2]$ & 5 & 5 \\ 
$[0 , 0]$ & 8 & 8 \\ 
$[1 , -2]$ & 5 & 5 \\  
$[2 , -4]$ & 1 & 1 \\ 
$[-3 , 3]$ & 2 & 2 \\ 
$[-2 , 1]$ & 5 & 5 \\ 
$[-1 , -1]$ & 5 & 3 \\ 
$[-4 , 2]$ & 1 & 1 \\ [1ex] 
\hline 
\end{tabular}
\end{table}

\begin{table}[H]
\caption{Two candidate $F_{110}$ planes have the following weight decomposition} 
\centering 
\begin{tabular}{c c c} 
\hline\hline 
Weight & Dimension of weight space in $A^* \otimes B^*$ & Dimension of weight space in $F_{110}$  \\ [0.5ex]
\hline 
$[2 , 2]$ & 1 & 1 \\ 
$[3 , 0]$ & 2 & 2 \\ 
$[4 , -2]$ & 1 & 1 \\  
$[0 , 3]$ & 2 & 2 \\ 
$[1 , 1]$ & 5 & 5 \\ 
$[2 , -1]$ & 5 & 5 \\ 
$[3 , -3]$ & 2 & 2 \\  
$[-2 , 4]$ & 1 & 1 \\ 
$[-1 , 2]$ & 5 & 5 \\ 
$[0 , 0]$ & 8 & 8 \\ 
$[1 , -2]$ & 5 & 5 \\  
$[2 , -4]$ & 1 & 1 \\ 
$[-3 , 3]$ & 2 & 2 \\ 
$[-2 , 1]$ & 5 & 5 \\ 
$[-1 , -1]$ & 5 & 4 \\  
$[-4 , 2]$ & 1 & 1 \\ 
$[-3 , 0]$ & 2 & 1 \\ [1ex] 
\hline 
\end{tabular}
\end{table}

The computation to produce these 5 candidate $F_{110}$ planes took extensive time in some cases, due to the parameters creating a difficult groebner basis computation when determining whether an $F_{110}$ plane passes $(210)$-test. The large number of candidates was not as much of a computational issue as all the $(210)$-tests can be parallelized. Some of these computations were done on Texas A{\&}M's High Performance Research Cluster as well as Texas A{\&}M's Math Department Cluster.

Using the skew-symmetry of $T_{\sln}$, we are able to produce $F_{011}$ and $F_{101}$ candidate weight spaces from the candidate $F_{110}$ spaces. A computer calculation verified that for each candidate triple $F_{110},F_{011},F_{101}$, the rank condition is not met for the $(111)$-test and consequently, there are no candidate $F_{111}$ spaces. Therefore, the rank of $T_{\slt}$ is greater than 15 and {this proves Theorem \ref{sl3_lower}, $\BR(T_{\slt}) \geq 16$.} This result is significant as it is the first example of an explicit tensor such that the border rank is at least $2m$ when $m<13$.

%
%
%
%

\subsection{Upper Bounds}

A numerical computer search has given a rank 20 decomposition of $T_{\slt}$. The technique used was a combination of Newton's Method and Lenstra–Lenstra–Lovász Algorithm to find rational approximations \cite{NewtonLLL}. {This technique formulated the problem as a nonlinear optimization problem that was solved to machine precision and then subsequently modified using the Lenstra-Lenstra-Lovász Algorithm to generate a precise solution with algebraic numbers given the numerical solution. As $T_{\slt} \in \C^8 \otimes \C^8 \otimes \C^8$, a rank 20 decomposition consists of finding $a_i,b_i,c_i \in \C^8$ such that $T = \sum_{i=1}^{20} a_i \otimes b_i \otimes c_i$. We take each vector $a_i,b_i,c_i$ to be a vector 8 variables, and using properties of elements of tensor products, we can multiply out the right hand side and have a system of equations for each entry of the tensor. This amounts to solving a system of 512 polynomial equations of degree 3 in 480 variables. We then use Newton's method to find roots to this system of equations. If it appears to converge to a solution, then we compute it to machine precision and use Lenstra-Lenstra-Lovász to find an algebraic solution that satisfies the initial polynomial conditions.} Therefore, {this decomposition proves Theorem \ref{sl3_upper}, $\TR(T_{\slt}) \leq  20$.}

Let $\zeta_6$ denote a primitive 6th root of unity. The following is the rank 20 decomposition of $T_{\slt}$. {One may verify that this is in fact a rank decomposition by showing that it satisfies the polynomial equations described above.} {We note that the below decomposition has some localy symmetry. For example, the first three terms share $ \begin{bmatrix}\Transpose{ 0 & 0 & 0 \\ \z& 0 & 1\\ 0 & 0 & 0 }\end{bmatrix}$ as a common factor with standard cyclic $\Z_3$-symmetry, i.e. this element occurs once in each tensor factor in the first three terms. We have grouped terms in the below presentation of the decomposition by this local symmetry. We were unable to determine any global symmetry in this presentation of this decomposition.}

\begin{align*}
&T_{\slt} = \\
&(\frac{1}{3^4 2}) \begin{bmatrix}\Transpose{ 0 & 0 & 0 \\ \z& 0 & 1\\ 0 & 0 & 0 }\end{bmatrix} \otimes \begin{bmatrix}\Transpose{ -6 & 9\z & 6\z^2 \\ -4\z^2 & 18 & 0 \\ -6\z & 0 & -12 }\end{bmatrix} \otimes \begin{bmatrix} 6\z & -4 & 0 \\ 9\z^2 & 0 & -9 \\ 0 & 4\z^2 & -6\z \end{bmatrix} + \\ 
&(\frac{1}{3^4 2^3}) \begin{bmatrix}\Transpose{ -6 & 9\z & 6\z^2 \\ -4\z^2& 18 &0\\ -6\z & 0 & -12} \end{bmatrix} \otimes \begin{bmatrix} \Transpose{0 & 0 & 0 \\ \z & 0 & 1 \\ 0 & 0 & 0 }\end{bmatrix} \otimes \begin{bmatrix} -6\z & 4 & -6\z^2 \\ 0 & 6\z & 9 \\ -6 & 0 & 0 \end{bmatrix} + \\ 
&(\frac{1}{3^4 2}) \begin{bmatrix} \Transpose{6 & 9\z & 0 \\ 4\z^2 & 0 & 4\z \\ 0 & 9\z^2 & -6} \end{bmatrix} \otimes \begin{bmatrix} \Transpose{12\z & 0 & 6 \\ -4 & -18\z & 0 \\ 6\z^2 & -9 & 6\z } \end{bmatrix} \otimes \begin{bmatrix} 0 & \z & 0 \\ 0 & 0 & 0 \\ 0 & 1 & 0 \end{bmatrix} +\\ 
&(\frac{\z}{3^5 2}) \begin{bmatrix}\Transpose{ 0 & -18 \z & -6 \\ 0& -18\z^2 & -4 \z \\ 6\z & -9 & 18\z^2 }\end{bmatrix} \otimes \begin{bmatrix}\Transpose{ 18\z^2 & 9\z & -6 \\ 0 & -18\z^2 & -4\z \\ 6 \z & 18 & 0 }\end{bmatrix} \otimes \begin{bmatrix} 0 & \z^2 & 0 \\ 0 & 0 & 0 \\ 0 & -1 & 0 \end{bmatrix} + \\ 
&(\frac{1}{ 3^4 2^2}) \begin{bmatrix}\Transpose{ -6\z & 9 & -6\z^2 \\ 4\z^2 & 18\z &0\\ -6 & 0 & -12\z} \end{bmatrix} \otimes \begin{bmatrix} \Transpose{0 & 0 & 0 \\ \z^2 & 0 & -1 \\ 0 & 0 & 0 }\end{bmatrix} \otimes \begin{bmatrix} -12\z & 0 & -6 \\ -9 & 6\z & 18\z^2 \\ -6\z^2 & 4 & 6\z \end{bmatrix} + \\ 
&(\frac{\z}{3^4}) \begin{bmatrix}\Transpose{0 & 0 & 0 \\ 1 & 0 & \z \\ 0 & 0 & 0 }\end{bmatrix} \otimes \begin{bmatrix} \Transpose{-3\z^2 & 0 & 1 \\ -2 & 9\z^2 & 0 \\ -3\z & 0 & -6\z^2 }\end{bmatrix} \otimes \begin{bmatrix} 6\z & 0 & -6 \\ 1 & 6\z & -9\z^2 \\ -6\z^2 & 4 & -12\z \end{bmatrix} + \\ 
&(\frac{1}{3^4 2^3}) \begin{bmatrix}\Transpose{ -12 & 0 & 6\z^2 \\ -4\z^2 & 18 &0\\ -6\z & -9\z^2 & -6 }\end{bmatrix} \otimes \begin{bmatrix} \Transpose{ -6\z^2 & 9 & 0 \\ 4\z & 0 & 4 \\ 0 & 9\z & 6\z^2 }\end{bmatrix} \otimes \begin{bmatrix} 0 & -4 & 6\z^2 \\ -9\z^2 & -6\z & 0 \\ 6 & 0 & 6\z \end{bmatrix} + \\ 
&(\frac{1}{ 3^5 2^3}) \begin{bmatrix} \Transpose{-18 & 9\z^2 & -6\z \\ 0& 18 & -4\z^2\\ 6\z^2 & 18\z & 0 }\end{bmatrix} \otimes \begin{bmatrix} \Transpose{0 & 18\z^2 & 6\z \\ 0 & -18 & 4\z^2 \\ -6\z^2 & 9\z & 18} \end{bmatrix} \otimes \begin{bmatrix} 0 & 4\z & 6\z^2 \\ -9\z^2 & 6 & 0 \\ -6\z & 0 & -6 \end{bmatrix} + \\
&(\frac{1}{3^3 2^3}) \begin{bmatrix}\Transpose{ -12\z & 0 & -6 \\ 8 & 18\z & -4\z^2\\ -6\z^2 & 9 & -6\z} \end{bmatrix} \otimes \begin{bmatrix}\Transpose{ -2\z^2 & 3 & 0 \\ 0 & 0 & 0 \\ 0 & 3\z & 2\z^2 }\end{bmatrix} \otimes \begin{bmatrix} 0 & 0 & -6\z \\ 9\z & 6 & 0 \\ 6\z^2 & 4\z & -6 \end{bmatrix} + \\ 
&(\frac{1}{3^5 2^3}) \begin{bmatrix}\Transpose{ 18\z^2 & 9\z & -6 \\ 12& -18\z^2 & 8\z \\ 6\z & 18 & 0 }\end{bmatrix} \otimes \begin{bmatrix}\Transpose{ 0 & 18\z & 6 \\ -6 & 18\z^2 & -2\z \\ -6\z & 9 & -18\z^2} \end{bmatrix} \otimes \begin{bmatrix} 0 & 0 & 6\z \\ -9\z & -6\z^2 & 0 \\ -6 & -4\z & 6\z^2 \end{bmatrix} + 
\end{align*}

\newpage

\begin{align*}
&(\frac{1}{3^4 2^3}) \begin{bmatrix} \Transpose{12 & 0 & -6 \\ 4& -18 &0\\ -6 & 9 & 6} \end{bmatrix} \otimes \begin{bmatrix} \Transpose{-6 & 9 & 0 \\ 4 & 0 & -4 \\ 0 & -9 & 6} \end{bmatrix} \otimes \begin{bmatrix} 0 & 0 & 6 \\ -9 & 6 & 0 \\ 6 & -4 & -6 \end{bmatrix} + \\ 
&(\frac{1}{3^3 2^3}) \begin{bmatrix}\Transpose{ -12 & 0 & 6 \\ 0& 18 & -4 \\ 6 & -9 & -6 }\end{bmatrix} \otimes \begin{bmatrix}\Transpose{ 2 & -3 & 0 \\ 0 & 0 & 0 \\ 0 & 3 & -2} \end{bmatrix} \otimes \begin{bmatrix} 0 & 4 & -6 \\ 9 & -6 & 0 \\ -6 & 0 & 6 \end{bmatrix} + \\ 
&(\frac{1}{3^3 2}) \begin{bmatrix} \Transpose{-2& 3 & 0 \\ 0& 0 &0\\ 0 & -3 & 2 }\end{bmatrix} \otimes \begin{bmatrix} \Transpose{12 & 0 & -6 \\ 0 & -18 & 4 \\ -6 & 9 & 6 }\end{bmatrix} \otimes \begin{bmatrix} 0 & -1 & 0 \\ 0 & 0 & 0 \\ 0 & 1 & 0 \end{bmatrix} + \\ 
&(\frac{1}{ 3^4 2}) \begin{bmatrix}\Transpose{ 6 & 9 & -6 \\ 4 & -18 &0\\ -6 & 0 & 12 }\end{bmatrix} \otimes \begin{bmatrix} \Transpose{0 & 0 & 0 \\ 1 & 0 & -1 \\ 0 & 0 & 0 }\end{bmatrix} \otimes \begin{bmatrix} 6 & 0 & -6 \\ 0 & -6 & 9 \\ -6 & 4 & 0 \end{bmatrix} + \\ 
&(\frac{1}{3^3 2}) \begin{bmatrix} \Transpose{0 & 0 & 0 \\  1 & 0 &-1\\ 0 & 0 & 0} \end{bmatrix} \otimes \begin{bmatrix} \Transpose{-6 & -9 & 6 \\ -4 & 18 & 0 \\ 6& 0 & -12} \end{bmatrix} \otimes \begin{bmatrix} 2 & 0 & 0 \\ -3 & 0 & 3 \\ 0 & 0 & -2 \end{bmatrix} +\\ 
&(\frac{1}{3^3 2^5}) \begin{bmatrix} \Transpose{-4 & 0 & 9 \\ -6 & 0 &4\\ -6 & 0 & 4} \end{bmatrix} \otimes \begin{bmatrix} \Transpose{-4 & 0 & -9 \\ 6 & 0 & 4\\ -6 & 0 & 4} \end{bmatrix} \otimes \begin{bmatrix} -4 & -6 & 6 \\ 0 & 0 & 0 \\ 9 & -4 & 4 \end{bmatrix} + \\ 
&(\frac{1}{3^32^5}) \begin{bmatrix}\Transpose{ 4 & 0 & 9 \\ -6 & 0 & 4\\ -6 & 0 & -4 }\end{bmatrix} \otimes \begin{bmatrix}\Transpose{ -4 & 0 & 9 \\ -6 & 0 & -4 \\ 6 & 0 & 4 }\end{bmatrix} \otimes \begin{bmatrix} 4 & -6 & 6 \\ 0 & 0 & 0 \\ 9 & -4 & -4 \end{bmatrix}+\\ 
&(\frac{1}{3^3 2^5}) \begin{bmatrix}\Transpose{ -4 & 0 & 9 \\ -6& 0 &-4\\ 6 & 0 & 4 }\end{bmatrix} \otimes \begin{bmatrix}\Transpose{ -4 & 0 & -9 \\ 6 & 0 & -4 \\ 6 & 0 & 4 }\end{bmatrix} \otimes \begin{bmatrix} 4 & 6 & 6 \\ 0 & 0 & 0 \\ -9 & -4 & -4 \end{bmatrix} + \\ 
&(\frac{1}{3^3 2^5}) \begin{bmatrix} \Transpose{-4 & 0 & -9 \\ 6& 0 & 4\\ -6 & 0 & 4} \end{bmatrix} \otimes \begin{bmatrix} \Transpose{-4 & 0 & 9 \\ -6 & 0 & 4 \\ -6 & 0 & 4 }\end{bmatrix} \otimes \begin{bmatrix} 4 & -6 & -6 \\ 0 & 0 & 0 \\ 9 & 4 & -4 \end{bmatrix} + \\ 
&(\frac{2}{3^2}) \begin{bmatrix}\Transpose{ -1 & 0 & 0 \\ 0& 0 &0\\ 0 & 0 & 1 }\end{bmatrix} \otimes \begin{bmatrix} \Transpose{0 & 0 & 0 \\ 0 & 0 & -2 \\ 3 & 0 & 0} \end{bmatrix} \otimes \begin{bmatrix} 0 & -2 & 0 \\ 0 & 0 & 0 \\ 3 & 0 & 0 \end{bmatrix} 
\end{align*}

{In an attempt to find a smaller rank decomposition, we found numerical evidence suggesting that $\BR(T_{\slt}) \leq 18$. The above method was unable to determine exact algebraic numbers for it to be an honest border rank decomposition. We include the approximate border rank decomposition, which was obtained as a numerical solution to machine precision using Newton's method, in Appendix A. This decomposition is satisfies the equation $T_{\slt} = \sum_{k=1}^{18} a_k(t) \otimes b_k(t) \otimes c_k(t) + O(t)$ to a maximum error in each entry of $3.88578058618805\, 10^{-16}$ ($\ell_0$ error). It also is satisfied with a sum of squares error of $1.85900227125328\, 10^{-15}$ ($\ell_2$ error), which is the square root of the sum of the squares of all errors in each entry.}

%
%

\appendix

\section{Appendix}

\begin{lstlisting}
## Approximate Border Rank 18 Decomposition for sl3. B is list of length 18,
## Each element of B is 3 lists, one for each tensor factor

t = var('t')
decomp_sl3_br18 = [
[[0.8582131228193816*t^4, 1.0*t^3, -0.656183917735616, 0, -0.6991867952664118*t^4, 0.48158077142326267*t, 0, 0],
[0.4276000554886944, 0, -0.18145759099190728*t^-4, -0.615151864463497*t, 0.3663065358105463, 0.41017982352375143*t^-3, 0.32919233217347255*t^4, -0.37760897344619454*t^3],
[-0.1349207590097993*t^-4, -0.6002859548136603*t^-3, 1.0, 0.5115159306456755*t^-5, -0.1595527810915205*t^-4, -0.6048099252034903*t^-1, -0.7277947872836206*t^-8, -0.40421722543545996*t^-7]],
[[-1.0*t^4, 0.35415463745033937*t^3, 0.0064263290832379735, 0, -0.36545526641405424*t^4, -0.07853056373003943*t, -0.032960816562400096*t^8, 0],
[0, 0, 1.0*t^-4, -0.2566348901375672*t, -0.22472886423613508, 0, -1.5211311526049585*t^4, 0],
[0.8999601061346961*t^-4, 0, -0.9891997951816988, 0.7757989654047442*t^-5, 0.6720573676204581*t^-4, 0, -0.7680275258871294*t^-8, 0.23098457888352464*t^-7]],
[[1.4652283519838352*t^4, 0, -0.49981037975335263, -0.1095789691297215*t^5, 0, -0.39151017186719006*t, -0.0035109361876165053*t^8, 0],
[-0.0370198277154002, -0.4868421872551675*t^-1, -0.6112939642144736*t^-4, 1.0*t, 0, -0.2400052666531245*t^-3, 0.3994050510000636*t^4, 0.38578320751551654*t^3],
[-0.976481792017985*t^-4, -0.6067214925499242*t^-3, 0, -0.723941360626237*t^-5, -1.0*t^-4, 0.5548674236521388*t^-1, 0.43144009347718587*t^-8, 0]],
[[0, 0, 0, 0.7384182957633163*t^5, 0, 0.644445145257651*t, 1.0*t^8, 0.14298165951600136*t^7],
[0.48917042278155703, -1.1676021066814524*t^-1, 0.0742778280974771*t^-4, 0, 0, 0.43490564244204044*t^-3, 0, -0.35107756870823037*t^3],
[0.35864016315209*t^-4, 0.08753789324397686*t^-3, 0.09745671407220928, 0.8747382868849779*t^-5, -0.08844676171127677*t^-4, 0, 1.0*t^-8, 0.26767660269283405*t^-7]],
[[0, 0.7154062494455652*t^3, 0.753065093225573, 0.25916815694138623*t^5, 1.0*t^4, -0.13328651389721624*t, -0.43754268152319187*t^8, 0.33712746822644685*t^7],
[0.31175249346291917, 0.4893435772426664*t^-1, -0.6750030230291015*t^-4, 0, 0, 0, -1.0*t^4, 0.2730200083141182*t^3],
[-0.48768995751368266*t^-4, 0.3666405285647919*t^-3, -0.509119354246841, -0.4700930934822484*t^-5, -0.22797281366964972*t^-4, -0.35924834770877284*t^-1, 0.21574925698409564*t^-8, -0.4457498667959218*t^-7]],
[[-0.3726117655275953*t^4, -0.6850149555795361*t^3, -0.47485190359116686, -0.3845493471094187*t^5, 1.3294054010695648*t^4, 0, 0, -0.360876333169775*t^7],
[-0.20585599800494148, 0.14163424899135532*t^-1, -0.07154541239160381*t^-4, 0, 0, 0.6161959625923827*t^-3, -1.0*t^4, 1.1951737645411018*t^3],
[0.32862886361063015*t^-4, -1.0*t^-3, 0.02421044619636145, -0.10037637090787854*t^-5, 0.8394716140987267*t^-4, 0.6969300342438073*t^-1, 0.33915176005821984*t^-8, -0.8710663402014699*t^-7]],
[[-0.26018060869681375*t^4, -0.5986510147484245*t^3, 0, 0, 0.34975430620458303*t^4, 0.6295753676920407*t, -0.008909285185907506*t^8, -1.0*t^7],
[0.1869896012120452, 0.21074450964709926*t^-1, 0.15931993306826575*t^-4, -0.5916688936748106*t, -1.0, 0.4704225477242227*t^-3, 0.4764298640490715*t^4, 0],
[-0.11381180152788646*t^-4, 0, -0.6321964884450085, 1.1046891281785145*t^-5, 0, 1.0*t^-1, -0.39603240282608254*t^-8, 0.9263232781144852*t^-7]],
[[0.5249458233325552*t^4, 0.48647993414509433*t^3, 0.288317913855889, -0.6006289011129895*t^5, -0.43339360309318065*t^4, -0.9021189701254521*t, -0.3580507322918174*t^8, 0],
[0, -0.6140778368802516*t^-1, 0, 0.3157042904634412*t, 1.0, 0.3354809468306432*t^-3, -0.4344836650252679*t^4, 0.3640820338419208*t^3],
[0, 0.021642017184031522*t^-3, 0, 1.0*t^-5, -0.16272749957295446*t^-4, 0, 0.3302270185683934*t^-8, -0.5277803864906128*t^-7]],
[[0, 0.5078943260216121*t^3, -0.8437068467257393, -0.16779079767696245*t^5, 0, -0.22773097773605452*t, 1.0*t^8, 0],
[-0.6816689949383756, 0, -0.5858467479048473*t^-4, 0.4936500799168931*t, -1.0, -0.06688550512654293*t^-3, 0, 0],
[0, 0.6059010244748405*t^-3, 0, 1.3353428165296193*t^-5, 0, 0, -0.6434299512786738*t^-8, -0.4609346348556645*t^-7]],
[[0.5342376006675218*t^4, -0.6086422770257539*t^3, 0, 0.4884368142723877*t^5, -0.718333221563019*t^4, 0, 0, -1.0*t^7],
[0.4628574417338272, 0.8309327397168714*t^-1, 0.21142383552160535*t^-4, -1.0560448029801677*t, 0, -1.0*t^-3, 0, 0.5907349063234486*t^3],
[-0.19068598535726494*t^-4, 1.251113151532601*t^-3, 0, 0.801440153916977*t^-5, -0.4621883155369931*t^-4, 0.649922008892636*t^-1, 0.4517624031930135*t^-8, 0.6177441885560955*t^-7]],
[[-0.46119950535589427*t^4, 0.7406376078020875*t^3, 0.29057187128303724, 0, 0.6143929403505343*t^4, 0.4568692437943906*t, 1.0*t^8, 0.36703376745263416*t^7],
[0, 1.080760392781858*t^-1, -0.14934226596262576*t^-4, 0, 0.463851046688496, -0.5286724400096557*t^-3, 0, -0.5122850384466106*t^3],
[0, -0.17675493191050737*t^-3, -0.7935243857498888, 1.0525678466492765*t^-5, 0, 0, 0.8232503249794846*t^-8, 1.0*t^-7]],
[[-0.1128613231857281*t^4, -1.004601481243395*t^3, -0.18790258087897918, 0.40149916700646043*t^5, -0.14415337265479422*t^4, 0, 0, 1.0*t^7],
[-0.6862576195287183, 0.5664746389899703*t^-1, 0.2908675685084361*t^-4, 1.3714715312276369*t, 1.0, -0.5713223321771269*t^-3, 0.5522456890804313*t^4, 0],
[-0.3006833454545763*t^-4, -1.5673937490149197*t^-3, 0, 0, 0.38713909105007405*t^-4, 0.9955939462610987*t^-1, 0.5677563588587456*t^-8, 0.5353436694572539*t^-7]],
[[-0.6245100623794343*t^4, 0, -1.0, 0.5065225347182659*t^5, 0.5317253204886139*t^4, 0, 0, -0.3719669835035259*t^7],
[0.5699237870531472, -0.2500784148082814*t^-1, 0.442945128286494*t^-4, 0.45638933734578824*t, 1.3973251849132609, 0, 0, 0],
[0, 0.5427787886626944*t^-3, 0.6987561651998788, 1.0*t^-5, 0, 0.20324504499293194*t^-1, -0.21496324065341224*t^-8, -0.35477523806474964*t^-7]],
[[0, -0.8732161724035096*t^3, -0.2040186481697669, -1.2010646422303124*t^5, -0.23508531984800027*t^4, 0.5654807700289866*t, 0, 1.0788666817176849*t^7],
[-1.0, 0.23352368487408065*t^-1, 0.2848223349469312*t^-4, 0.1768073101741446*t, 1.1645543351076657, 0, 0, 0],
[0.4548846200775119*t^-4, 0, -1.0, 0.7165423092004359*t^-5, -0.2616313451578387*t^-4, -1.2932734528263643*t^-1, -0.30361754210884495*t^-8, 0.12055753515645255*t^-7]],
[[-0.36174656010597567*t^4, 0, 0, 1.0*t^5, 0, -0.9529930944532633*t, 0, 0],
[-0.249684891380413, -0.48049270284391143*t^-1, 0.27222663948587644*t^-4, 0.06402554598955673*t, 0.32897469345165686, -0.29230508678584777*t^-3, 0.3854984554913255*t^4, 1.084386049943073*t^3],
[0.5943543096084136*t^-4, 0, -0.4586866892423329, -1.0*t^-5, -0.2869864460438895*t^-4, 0.7692509696955846*t^-1, -0.5902387210012353*t^-8, -0.10600670730550378*t^-7]],
[[0.2878925845174677*t^4, 0.5281046432775963*t^3, -0.2234313203594912, 0, -0.14791845243016916*t^4, -0.33262458187432487*t, 1.108538091915835*t^8, 0.2489733522919065*t^7],
[-0.19980163401634468, -0.4519313364099573*t^-1, -0.9400278789365702*t^-4, 1.7023099438662026*t, 0.47743291480661737, 0.40637104558436565*t^-3, -0.9678003132293527*t^4, -1.0*t^3],
[1.0*t^-4, 0.391031885184312*t^-3, 0, -1.0*t^-5, 0.8654063566747124*t^-4, 0, 0, 0]],
[[1.0110775513705723*t^4, -0.3307689454870473*t^3, -0.45870287159609613, 0, -0.2918463909613691*t^4, 0.30167345561087666*t, 0.582399773106054*t^8, 0],
[0, -1.0*t^-1, 0.7191039431350641*t^-4, 0.5392536432576923*t, -0.39197591630294476, -0.7143544436146056*t^-3, 0, 0.8086746481540814*t^3],
[0.4561800893644289*t^-4, -0.7356641704491395*t^-3, 0.4710082951189216, 0, 0, 0.7109669616748934*t^-1, -1.0*t^-8, 0]],
[[-0.41687446856120663*t^4, 0.3197337020122604*t^3, 0.3753986532686001, 0, -0.40495125990774045*t^4, 0.2979428965247865*t, -0.6344817450528121*t^8, -0.212014762523078*t^7],
[-0.2396929297153146, 0.7894128907681753*t^-1, 0, 0, 0.2995112979592937, 1.0*t^-3, 0.09278089610904368*t^4, 0.5160805291028968*t^3],
[0.4011063039077983*t^-4, 0, 0, 0.2503035377295029*t^-5, 1.0*t^-4, -0.9730216219558445*t^-1, 0.6637049194937615*t^-8, -0.5073657590422873*t^-7]]
]
\end{lstlisting}

\bibliography{ref}
\bibliographystyle{halpha}

\end{document}